\documentclass[12pt,reqno]{amsart}
\usepackage{enumerate,color,amsmath,amssymb,amsfonts,bm}
\usepackage{frcursive,fullpage,graphicx,psfrag}
\usepackage{mathrsfs}

\newtheorem{Remark}{Remark}
\newtheorem{Theorem}{Theorem}
\newtheorem{Proposition}{Proposition}
\newtheorem{Corollary}{Corollary}
\newtheorem{Conjecture}{Conjecture}

\def\bbone{{\mathchoice {\rm 1\mskip-4mu l} {\rm 1\mskip-4mu l}
{\rm 1\mskip-4.5mu l} {\rm 1\mskip-5mu l}}}
\newcommand{\ux}{\mathbf{x}}
\newcommand{\ups}{\Upsilon}

\begin{document}

\author{Abdelmalek Abdesselam}
\address{Abdelmalek Abdesselam, Department of Mathematics,
P. O. Box 400137,
University of Virginia,
Charlottesville, VA 22904-4137, USA}
\email{malek@virginia.edu}

\title{An algebraic independence result related to a conjecture of Dixmier on binary form invariants}

\begin{abstract}
In order to better understand the structure of classical
rings of invariants for binary forms, Dixmier proposed, as a conjectural
homogeneous system of parameters,
an explicit collection of invariants previously studied by Hilbert.
We generalize Dixmier's collection and show that a particular subfamily is algebraically independent. Our proof relies on showing certain alternating sums of products of binomial coefficients are nonzero. Along the way we provide a very elementary proof {\it \`a la} Racah, namely, only using the Chu-Vandermonde Theorem, for Dixon's Summation Theorem. 
We also provide explicit computations of invariants, for the binary octavic, which can serve as ideal introductory examples to Gordan's 1868 method in classical invariant theory.
\end{abstract}

\maketitle

\tableofcontents

\section{Introduction}

Throughout this article we will work over the field $\mathbb{C}$ of complex numbers. 
For an integer $d\ge 1$, we will denote by $S_d$ the $(d+1)$-dimensional vector space of binary forms
\[
F(\ux)=\sum_{i=0}^{d}\ f_i\ x_1^{d-i} x_2^{i}\ ,
\]
namely, homogeneous polynomials of degree $d$ in the pair of variables $\ux=(x_1,x_2)$.
This space carries a natural left $SL_2$ action defined as follows.
For a matrix
\[
g=\left(
\begin{array}{cc}
g_{11} & g_{12} \\
g_{21} & g_{22}
\end{array}
\right)
\]
in $SL_2$, we let 
$g\cdot\ux=(g_{11}x_1+g_{12}x_2,g_{21}x_1+g_{22}x_2)$
and
$(g\cdot F)(\ux)=F(g^{-1}\cdot \ux)$.
A polynomial $C(F,\ux)=C(f_0,\ldots,f_d,x_1,x_2)$ is called a covariant of the generic
binary form of degree $d$ if it identically satisfies
$C(g\cdot F,g\cdot \ux)=C(F,\ux)$
for all $g\in SL_2$.
Such covariants form a ring
${\rm Cov}_d\subset \mathbb{C}[f_0,\ldots,f_d,x_1,x_2]$
which is bigraded by $({\rm degree},{\rm order})$ where ``degree'' refers to the degree in the $f$ variables and ``order'' means the degree in the $x$ variables.
The order zero subring is the ring of invariants ${\rm Inv}_d$.
The study of these rings is a classical subject in mathematics.

Minimal systems of bihomogeneous generators for ${\rm Cov}_5$ (23 generators)
and ${\rm Cov}_6$ (26 generators) were obtained by Gordan in~\cite{Gordan}.
A system of generators for ${\rm Cov}_8$ was obtained by von Gall
in~\cite{vonGall1,vonGall2,vonGall3}. 
He also treated the more difficult case of ${\rm Cov}_7$ in~\cite{vonGall4}.
Note that when a generating system of ${\rm Cov}_d$ is known, this
immediately gives one for ${\rm Inv}_d$ by keeping generators of order zero.
Moreover, in the classical approach of Gordan and von Gall, the study of ${\rm Inv}_d$
necessitates that of the entire ring ${\rm Cov}_d$. 
Progress in the study of these rings has stagnated for a long time with a few notable exceptions.
Shioda rederived a minimal system for ${\rm Inv}_8$ and also found the syzygies among these generators~\cite{Shioda}.
von Gall's system for ${\rm Cov}_7$ was generating but not minimal. Six elements in his list were in fact reducible.
The determination of a truly minimal system of 147 generators for ${\rm Cov}_7$ has been done in~\cite[p. 227]{Croni} and~\cite{Bedratyuk}.
Previously, some doubt left by von Gall's work about the number of generators for
${\rm Inv}_7$ was remedied by Dixmier and Lazard~\cite{DixmierL}. 
Only very recently, Brouwer and Popoviciu obtained the minimal systems of generators for
${\rm Inv}_9$ (92 invariants) in~\cite{BrouwerP1} and for  ${\rm Inv}_{10}$ 
(106 invariants) in~\cite{BrouwerP2}.
Finally, and even more recently, in a tour de force of computational algebra
Lercier and Olive~\cite{Olive,LercierO} managed to go beyond von Gall's 1888 results for covariants and determined minimal
systems of generators for ${\rm Cov}_9$
(476 covariants) and ${\rm Cov}_{10}$ (510 covariants).
These rings are notoriously difficult to present explicitly. The evident complexity displayed by the above examples is supplemented by general results of Kac~\cite{Kac}
and Popov~\cite{Popov} from which one should expect high numbers of generators as well as a high homological dimension.
Yet and despite this complexity, the study of rings of invariants and covariants of binary forms still is a useful task, partly because of connections to rings of Siegel modular forms (see~\cite{Igusa} for ${\rm Inv}_d$ and the more recent~\cite{CleryFvdG}
for ${\rm Cov}_d$).
In an effort to find some regularity in the chaotic structure of the rings ${\rm Inv}_d$, Dixmier~\cite{Dixmier}
proposed three conjectures about their Hilbert series, according to the congruence of $d$ mod 4. In the particular case where $d$ is divisible by 4, he conjectured an explicit homogeneous system of parameters (HSOP) with degree sequence $(2,3,\ldots,d-1)$. A recent study of such HSOP's for ${\rm Inv}_d$ can also be found in~\cite{BrouwerDP}. Another investigation related to Dixmier's work is~\cite{Derksen}.

For the reader's convenience, we now recall the definitions and facts from commutative algebra about HSOP's, regular sequences, etc. Our references for this material are: Brion's lectures~\cite[\S3.4]{Brion}, the book by Bruns and
Herzog~\cite[Ch. 6]{BrunsH} and the one by Derksen and
Kemper~\cite[\S2.4 and \S2.5]{DerksenK}.
Let $R=\oplus_{j\ge 0} R_j$ be a graded $\mathbb{C}$-algebra with $R_0=\mathbb{C}$, and $\theta_1,\ldots,\theta_q$ be a sequence of homogeneous elements in $R$. This sequence is called a HSOP if $\theta_1,\ldots,\theta_q$ are algebraically independent over $\mathbb{C}$
and $R$ is a finitely generated module over the subring
$\mathbb{C}[\theta_1,\ldots,\theta_q]$.
Assume $R$ is finitely generated as a $\mathbb{C}$-algebra, then thanks to Hilbert's 
homogeneous version of Noether's normalization lemma, HSOP's are guaranteed to exist.
The length $q$ of such a HSOP must coincide with the Krull dimension of $R$.
If $R$ is a {\em free} module over $\mathbb{C}[\theta_1,\ldots,\theta_q]$  for some HSOP
$\theta_1,\ldots,\theta_q$ then $R$ is called a Cohen-Macaulay algebra. This free module property will then hold for {\em any} HSOP. Moreover, in this case a HSOP
must also be a regular sequence in $R$.
Recall that an $R$-regular sequence or simply a {\em regular sequence in} $R$ is a sequence of elements $\theta_1,\ldots,\theta_s$
in $R$ such that $R/\langle\theta_1,\ldots,\theta_s\rangle\neq 0$ and for all $i=1,\ldots,s$, multiplication by $\theta_i$ in $R/\langle\theta_1,\ldots,\theta_{i-1}\rangle$
is injective. For $R={\rm Inv}_d$, $q={\rm dim}(S_d)-{\rm dim}(SL_2)=(d+1)-3=d-2$, since the stabilizer of a generic $F\in S_d$ is finite.
By the Hochster-Roberts Theorem, $R={\rm Inv}_d$ is Cohen-Macaulay. Moreover, a regular sequence $\theta_1,\ldots,\theta_s$ in ${\rm Inv}_d$ must have $s\le q=d-2$ and can always be completed into a HSOP. One should also mention a criterion due to Hilbert which characterizes sequences of homogeneous elements $\theta_1,\ldots,\theta_{d-2}$
that are HSOP's. Namely, to be an HSOP, a necessary and sufficient condition is that $\theta_1,\ldots,\theta_{d-2}$ provide set-theoretic equations for the nullcone
$\mathscr{N}_d=\{F\in S_d\ |\ \forall J\in{\rm Inv}_d, J(F)=0\}$.
Finally, since this notion will play a role later in this article, let us call two sequences of homogeneous elements $\theta_1,\ldots,\theta_s$ and $\theta'_1,\ldots,\theta'_s$ {\em triangularly related} if for all $i$ one can write $\theta'_i=\gamma_i\theta_i+P_i(\theta_1,\ldots,\theta_{i-1})$ for some nonzero scalars $\gamma_i\in\mathbb{C}$ and for some polynomials $P_i$ with complex coefficients. Clearly, this is an equivalence relation which preserves the HSOP and regular sequence properties.

Recall that for $0\le k\le \min(m,n)$ and for $F\in S_m$ and $G\in S_n$, one has the classical notion of transvectant $(F,G)_k\in S_{m+n-2k}$ given by
\begin{equation}
(F,G)_k=\frac{(m-k)!(n-k)!}{m!\ n!}
\ \sum_{l=0}^{k}(-1)^l\ \left(\begin{array}{c}k\\ l\end{array}\right)
\frac{\partial^k F}{\partial x_1^{k-l}\partial x_2^{l}}\ 
\frac{\partial^k G}{\partial x_1^{l}\partial x_2^{k-l}}\ .
\label{transdef}
\end{equation}
This bilinear operation realizes the $SL_2$-equivariant projection
$S_m\otimes S_{n}\rightarrow S_{m+n-2k}$ of the Clebsch-Gordan decomposition.
We now assume $d=2k$ with $k$ even and take $n\ge k$. For fixed $F\in S_d$, consider the linear map
\[
\begin{array}{llll}
\mathcal{L}_{n}^{F}: & S_n & \longrightarrow & S_n \\
 & G & \longmapsto & (F,G)_k
\end{array}
\]
and for $1\le p\le n+1$ we
denote by $\mathscr{H}_{n,p}(F)$ the coefficient of $\lambda^p$ in the characteristic polynomial
${\rm det}(\lambda {\rm Id}-\mathcal{L}_{n}^{F})$ of the map $\mathcal{L}_{n}^{F}$.
These polynomials in $F$ were introduced by Hilbert in his
K\"onigsberg Habilitationsschrift~\cite{Hilbert} where he showed they are $SL_2$ invariants of $F$ (see also~\cite[\S3]{LittelmannP}).

\begin{Conjecture} (Dixmier~\cite[Conjecture 3']{Dixmier})\ 

When $d$ is divisible by 4,
$\mathscr{H}_{d-2,2}(F), \mathscr{H}_{d-2,3}(F),\ldots, \mathscr{H}_{d-2,d-1}(F)$
form a HSOP for ${\rm Inv}_d$.
\end{Conjecture}

An equivalent reformulation of Dixmier's Conjecture is that, for $d$ divisible by 4,
\[
\mathscr{P}_{d-2,2}(F), \mathscr{P}_{d-2,3}(F),\ldots, \mathscr{P}_{d-2,d-1}(F)
\]
form a HSOP for ${\rm Inv}_d$, where
\[
\mathscr{P}_{n,p}(F)={\rm tr}[(\mathcal{L}_{n}^{F})^p]\ .
\]
Indeed, $\mathscr{H}_{d-2,2}(F), \mathscr{H}_{d-2,3}(F),\ldots, \mathscr{H}_{d-2,d-1}(F)$ and $\mathscr{P}_{d-2,2}(F), \mathscr{P}_{d-2,3}(F),\ldots, 
\mathscr{P}_{d-2,d-1}(F)$ are triangularly related. This follows from the classical explicit formulas relating the elementary symmetric functions $\pm\mathscr{H}_{n,p}$
of the eigenvalues of $\mathcal{L}_{n}^{F}$ to the power sums
$\mathscr{P}_{n,p}$.
Note that one has identically
$\mathscr{H}_{n,1}=\mathscr{P}_{n,1}=0$ since binary forms have no linear invariant.
As mentioned above, a HSOP for
${\rm Inv}_{d}$ must have exactly $d-2$ elements. This explains why Dixmier picked the sequence corresponding to $n=d-2$ for his conjecture. We propose to enlarge Dixmier conjecture as follows.

\begin{Conjecture}\label{genconj}
For $d=2k$ with $k$ even and for $k\le n\le 2k-2$,
$\mathscr{P}_{n,2}, \mathscr{P}_{n,3},\ldots,\mathscr{P}_{n,n+1}$
is a regular sequence in the ring ${\rm Inv}_{d}$.
\end{Conjecture}

As will be made clear in \S5, Dixmier's Conjecture is very difficult. Our Conjecture
\ref{genconj} suggests a more progressive approach of establishing the regular sequence property for increasing values of $n$.
The main result of this article is the following modest step in this direction.

\begin{Theorem}\label{mainthm}
With the same hypotheses and notation as in Conjecture~\ref{genconj},
we have that
$\mathscr{P}_{k,2}, \mathscr{P}_{k,3},\ldots,\mathscr{P}_{k,k+1}$
are algebraically independent.
\end{Theorem}

The main body of the proof of this theorem is in \S2 where we compute the Jacobian matrix of the invariants at a suitable point or binary form $\mathbb{F}$. Algebraic independence follows from this matrix having full rank which amounts to showing the nonvanishing of some combinatorial sums $\ups_m$. This nonvanishing is established in \S3 by revisiting recent work of Guo, Jouhet and Zeng~\cite{GuoJZ}. In fact, one can write $\ups_m$ as an explicit sum of positive terms~\cite[theorem 1.2]{GuoJZ} which itself is a consequence of a rather intimidating $q$-hypergeometric
multisum identity of Andrews~\cite[Theorem 4]{Andrews}.  
In \S3, we follow a more elementary approach via a recursive formula~\cite[Lemma 2.1]{GuoJZ}.
This gives us an opportunity to introduce several improvements on the derivation
in~\cite{GuoJZ}: we do not use the Pfaff-Saalsch\"utz identity but only the simpler
Chu-Vandermonde Theorem. We start the recursion at $m=2$ instead of $m=3$ which forced the
authors of~\cite{GuoJZ} to invoke rather than deduce Dixon's Theorem. As a pleasant surprise, we obtained a proof of the latter which is very elementary and perhaps new.
Note that the computations in \S2 were originally done using the graphical calculus developed
in~\cite[\S2]{Abdesselam12}. The derivatives with respect to $f_s$, $0\le s\le d$,
of an invariant $\mathscr{P}$ are best packaged into its first evectant (see~\cite[\S5]{Chipalkatti} for a definition). We then wrote down a graphical formula for the covariant $\mathscr{W}$
obtained as the homogeneous Wronskian of these evectants (see~\cite{AbdesselamC7}
for a recent study of such Wronskians from an invariant-theoretic perspective).
The existence of a point where the Jacobian matrix has full rank is equivalent to 
the covariant $\mathscr{W}$ not being identically zero. With the help of our graphical representation for $\mathscr{W}$, and also with some inspiration taken
from~\cite[\S3]{LittelmannP}, 
we found a suitable point of specialization $\mathbb{F}$ in the nullcone.
However, in order to make our proof accessible to a wider audience,
we erased our footsteps in our writing of \S2 which can be read without knowledge
of~\cite[\S2]{Abdesselam12} and which
only requires elementary linear algebra and multivariable calculus.
In \S4, we explicitly compute the invariants $\mathscr{P}$ that are relevant for
Theorem \ref{mainthm}, in the case of the binary quartic (trivial) and that of the binary
octavic (involved but instructive).
In \S5,
which assumes some familiarity with the graphical calculus of~\cite[\S2]{Abdesselam12},
we provide a hopefully insightful discussion of Dixmier's Conjecture by extracting
some of the combinatorial difficulties it contains.

\section{Reduction to a nonvanishing statement for some combinatorial sums}

Throughout this article, we will, similarly to Iverson's bracket, use the notation
$\bbone\{\cdots\}$ for the indicator function of the condition between braces.
Let $\mathcal{B}$ denote the basis of monomials $x_1^{k-i}x_2^{i}$, $0\le i\le k$ for the space $S_k$. For a linear operator $M:S_k\rightarrow S_k$ we will denote by $[M]_{ij}$ the matrix elements of this operator in the basis $\mathcal{B}$.
Let $\mathscr{J}(F)=
(\mathscr{J}_{r,s})_{\substack{{2\le r\le k+1}\\{0\le s\le 2k}}}$
denote the Jacobian matrix of the invariants $\mathscr{P}_{k,2},\ldots,
\mathscr{P}_{k,k+1}$ at some binary form $F\in\ S_{d}$.
Namely,
\[
\mathscr{J}_{r,s}=\frac{\partial}{\partial f_s}
\mathscr{P}_{k,r}(F)\ .
\]
For $0\le s\le d=2k$ and $r\ge 1$, by the multivariate chain rule and the cyclic property of the trace we have
\begin{equation}
\frac{\partial}{\partial f_s}\mathscr{P}_{k,r}(F)=
r\times\ \sum_{i,j=0}^{k}\ [(\mathcal{L}_{k}^{F})^{r-1}]_{ij}
\ \frac{\partial}{\partial f_s}[\mathcal{L}_{k}^{F}]_{ji}\ .
\label{chainrule}
\end{equation}
A straightforward computation using (\ref{transdef})
gives, in general for $0\le i\le m$, $0\le j\le n$ and $0\le k\le \min(m,n)$,
\[
(x_1^{m-i}x_2^{i},x_1^{n-j}x_2^{j})_k=T_{i,j}^{m,n,k}\ x_1^{m+n-k-i-j} x_2^{i+j-k}
\]
where
\begin{eqnarray*}
T_{i,j}^{m,n,k} & = & \frac{(m-k)!\ (n-k)!\ k!\ i!\ (m-i)!\ j!\ (n-j)!}{m!\ n!} \\
 & & \times
\sum_{l=\max(0,k-j,k-m+i)}^{\min(i,n-j,k)}
\frac{(-1)^l}{l!(k-l)!(i-l)!(n-j-l)!(j-k+l)!(m-i-k+l)!}\ .
\end{eqnarray*}
In the application of this formula to our case of interest, the sum reduces to a single term and therefore, when $0\le i,j\le k$,
\[
[\mathcal{L}_{k}^{F}]_{ij}=f_{i-j+k}\times T_{i-j+k,j}^{2k,k,k}
\]
with
\[
T_{i-j+k,j}^{2k,k,k}=(-1)^{j}\times\frac{k!\ (k-i+j)!\ (k+i-j)!}{(2k)!\ i!\ (k-i)!}
\]
since $k$ is even.
We now specialize our calculation of the Jacobian matrix to the particular unstable form
$\mathbb{F}=x_1^{k-1}x_2^{k+1}$
which has as coefficients
$f_s=\ \bbone\{s=k+1\}$,
for $0\le s\le 2k$.
Thus
\[
[\mathcal{L}_k^{\mathbb{F}}]_{ij}=\ \bbone\{i=j+1\}\times\ 
\frac{(-1)^{j}\binom{k}{j+1}}{\binom{2k}{k+1}}
\]
and the corresponding matrix essentially is a nilpotent Jordan-like matrix with nonzero entries only immediately below the diagonal.
One also has
\[
\left.\frac{\partial}{\partial f_s}[\mathcal{L}_{k}^{F}]_{ji}\right|_{F=\mathbb{F}}
=\ \bbone\{s=j-i+k\}\times T_{j-i+k,i}^{2k,k,k}\ .
\]
Let $\mathbb{J}=
(\mathbb{J}_{r,s})_{\substack{{2\le r\le k+1}\\{0\le s\le 2k}}}$
denote the Jacobian matrix $\mathscr{J}(F)$ specialized at $F=\mathbb{F}$.
Then an immediate computation using (\ref{chainrule})
gives
\[
\mathbb{J}_{r,s}=
\ \frac{r\times(-1)^{\binom{r}{2}}\times \bbone\{s=k-r+1\}}
{\binom{2k}{k+r-1}
{\binom{2k}{k+1}}^{r-1}}\times
\mathscr{N}_{k,r}
\]
where
\begin{equation}
\mathscr{N}_{k,r}=
\sum_{j=0}^{k-r+1}(-1)^{rj}
\left(\begin{array}{c} k \\ j \end{array}\right)
\left(\begin{array}{c} k \\ j+1 \end{array}\right)\cdots
\left(\begin{array}{c} k \\ j+r-1 \end{array}\right)\ .
\label{Ndefeq}
\end{equation}
Consider the maximal $k\times k$ minor
\[
\widehat{\mathbb{J}}={\rm det}\left[
(\mathbb{J}_{r,s})_{\substack{{2\le r\le k+1}\\{0\le s\le k-1}}}\right]\ .
\]
Then $\widehat{\mathbb{J}}$ is equal to an obviously nonzero number times the product
$\prod_{r=2}^{k+1}\mathscr{N}_{k,r}$.

The algebraic independence in Theorem \ref{mainthm}
follows from $\widehat{\mathbb{J}}\neq 0$ which itself reduces to showing that
for all $r$, $2\le r\le k+1$, we have $\mathscr{N}_{k,r}\neq 0$.
The case where $r$ is even is of course trivial, whereas the more involved alternating sum
situation where $r$ is odd will be taken care of in Proposition \ref{nonvan} below.

\section{Some combinatorial identities including those of von Szily and Dixon}

In this section we will adopt the following convention regarding ratios of products of factorials. The $a$'s and $b$'s being elements of $\mathbb{Z}$, we let
\[
\frac{a_1!\cdots a_m!}{b_1!\cdots b_n!}:=a_1!\times\cdots\times a_m!\times
\frac{1}{b_1!}\times\cdots\times\frac{1}{b_n!}
\]
where, for any $n\in\mathbb{Z}$, we set by definition
\[
n!=\left\{
\begin{array}{cl}
{\rm ordinary}\ n! & {\rm if}\ n\ge 0\ ,\\
0 & {\rm if}\ n<0\ ,
\end{array}
\right.
\]
and
\[
\frac{1}{n!}=\left\{
\begin{array}{cl}
{\rm ordinary}\ \frac{1}{n!} & {\rm if}\ n\ge 0\ ,\\
0 & {\rm if}\ n<0\ .
\end{array}
\right.
\]
Beware that with such a convention $\frac{n!}{n!}$ is not equal to 1 but rather
$\bbone\{n\ge 0\}$.
For any $n,k\in\mathbb{Z}$, we also define the binomial coefficients
$\binom{n}{k}=\frac{n!}{k!(n-k)!}$ with the previous convention enforced. In particular, the coefficients are zero unless $0\le k\le n$.
In all the following combinatorial sums, the range of summation will be $\mathbb{Z}$ and therefore omitted. 

\begin{Proposition}(Chu-Vandermonde)\label{ChuV}
$\forall a,b,c\in\mathbb{Z}$,
\[
\sum_k \left(\begin{array}{c} a \\ k \end{array}\right)
\left(\begin{array}{c} b \\ c-k \end{array}\right)=
\bbone\left\{
\begin{array}{c}
a\ge 0\\
b\ge 0
\end{array}
\right\}
\left(\begin{array}{c} a+b \\ c \end{array}\right)\ .
\]
\end{Proposition}

\noindent{\bf Proof:}
Both sides of the equation vanish unless $a\ge 0$ and $b\ge 0$ which we now assume.
Taking the sum of the inequalities $k\ge 0$ and $c-k\ge 0$ implies $c\ge 0$.
So when $c<0$ both sides vanish. Taking the sum of the inequalities $k\le a$ and $c-k\le b$
implies $c\le a+b$.
So also when $c>a+b$ both sides vanish.
The remaining case is when $0\le c\le a+b$. Take $A$ a set of cardinality $a$. Take $B$ a set of cardinality $b$ that is disjoint from $A$. Finally, count subsets $C\subset A\cup B$ with cardinality $c$ by conditioning on the cardinality $k$ of $C\cap A$.
\qed

For our purposes a more practical form of the Chu-Vandermonde convolution identity is the following.

\begin{Corollary}\label{practicalCV}

$\forall a,b,c\in\mathbb{Z}$,
\[
\left(\begin{array}{c} 2a \\ a+c \end{array}\right)
\left(\begin{array}{c} 2b \\ b+c \end{array}\right)=
\frac{(2a)!(2b)!}{(a+b)!}\sum_l
\frac{1}{(a-l)!(b-l)!(l+c)!(l-c)!}\ .
\]
\end{Corollary}

\noindent{\bf Proof:}
The summand in the RHS includes the indicator functions
of $a-l\ge 0$ and $l-c\ge 0$ and therefore the sum of these inequalities gives
$a-c\ge 0$.
Likewise, the conditions $b-l\ge 0$ and $c+l\ge 0$ are included and so is their consequence
$b+c\ge 0$. 
We can thus multiply and divide by $(a-c)!$ and $(b+c)!$ which gives
\[
{\rm RHS}=\frac{(2a)!(2b)!}{(a+b)!(a-c)!(b+c)!}\times \sum_l
\left(\begin{array}{c} a-c \\ l-c \end{array}\right)
\left(\begin{array}{c} b+c \\ b-l \end{array}\right)\ .
\]
Changing variables to $k=l-c$ and applying Proposition \ref{ChuV}
gives
\[
{\rm RHS}=\frac{(2a)!(2b)!}{(a+b)!(a-c)!(b+c)!}\times
\bbone\left\{
\begin{array}{c}
a-c\ge 0\\
b+c\ge 0
\end{array}
\right\}
\left(\begin{array}{c} a+b \\ b-c \end{array}\right)
\]
which reduces to LHS after cleaning up the expression, taking our convention into account.
\qed

For all $m\ge 1$ and $a_1,\ldots,a_m\in\mathbb{Z}$, let
\[
\ups_m(a_1,\ldots,a_m)=\sum_k (-1)^k
\left(\begin{array}{c} 2a_1 \\ a_1+k \end{array}\right)
\left(\begin{array}{c} 2a_2 \\ a_2+k \end{array}\right)\cdots
\left(\begin{array}{c} 2a_m \\ a_m+k \end{array}\right)\ .
\]
These alternating combinatorial sums satisfy the following recursion due to
Guo, Jouhet and Zheng~\cite[Lemma 2.1]{GuoJZ}.

\begin{Proposition}

$\forall m\ge 2$, $\forall a_1,\ldots,a_m\in\mathbb{Z}$,
\begin{equation}
\ups_m(a_1,\ldots,a_m)=\sum_l
\frac{(2a_1)!\ (2a_2)!}{(a_1+a_2)!(2l)!(a_1-l)!(a_2-l)!}
\times \ups_{m-1}(l,a_3,\ldots,a_m)\ .
\label{recureq}
\end{equation}
\end{Proposition}

\noindent{\bf Proof:}
Using $Q$ as a placeholder for the product $\prod_{2<i\le m}\binom{2a_i}{a_i+k}$, we have
\[
\ups_m(a_1,\ldots,a_m)=\sum_k (-1)^k
\left(\begin{array}{c} 2a_1 \\ a_1+k \end{array}\right)
\left(\begin{array}{c} 2a_2 \\ a_2+k \end{array}\right)\times Q
\]
in which we insert the identity of Corollary \ref{practicalCV}
for the first two binomial coefficients
with the result
\[
\ups_m(a_1,\ldots,a_m)=\sum_k\sum_l
\frac{(-1)^k\ (2a_1)!\ (2a_2)!}{(a_1+a_2)!(a_1-l)!(a_2-l)!(l-k)!(l+k)!}\times Q\ .
\]
The summand includes the indicator function of the conditions $l-k\ge 0$ and $l+k\ge 0$
which imply $2l\ge 0$. It is thus legitimate to multiply and divide by $(2l)!$.
Also note that only finitely many pairs $(k,l)$ contribute to the last sum because of the implied conditions $0\le l\le \min(a_1,a_2)$ and $-l\le k\le l$.
Hence Fubini's Theorem applies and one can write
\[
\ups_m(a_1,\ldots,a_m)=\sum_l\sum_k
\frac{(-1)^k\ (2a_1)!\ (2a_2)!}{(a_1+a_2)!(a_1-l)!(a_2-l)!(2l)!}\times
\left(\begin{array}{c} 2l \\ l+k \end{array}\right)
\times Q\ .
\]
The result then follows from the definition of $\ups_{m-1}(l,a_3,\ldots,a_m)$.
\qed

\medskip
We now look at the first three simplest cases.

\begin{Proposition}

$\forall a_1\in\mathbb{Z}$, we have
\begin{equation}
\ups_1(a_1)=\bbone\{a_1=0\}\ .
\label{up1eq}
\end{equation}
\end{Proposition}

\noindent{\bf Proof:}
By changing variables to $j=a_1+k$, we have from the definition
\[
\ups_1(a_1)=(-1)^{a_1}\sum_j (-1)^j
\left(\begin{array}{c} 2a_1 \\ j \end{array}\right)\ .
\]
If $a_1<0$, the RHS is zero by convention. If $a_1=0$, the RHS
reduces to $\left(\begin{array}{c} 0 \\ 0 \end{array}\right)=1$. If $a_1>0$, Newton's
Binomial Theorem gives ${\rm RHS}=(-1)^{a_1}(1-1)^{a_1}=0$.
\qed

The case $m=2$ is already more interesting since it amounts to
the von Szily identity~\cite{vonSzily} for super-Catalan numbers (see, e.g.,~\cite{Gessel}).

\begin{Proposition}(von Szily)
$\forall a_1,a_2\in\mathbb{Z}$,
\begin{equation}
\ups_2(a_1,a_2)=\frac{(2a_1)!(2a_2)!}{(a_1+a_2)!a_1! a_2!}\ .
\label{up2eq}
\end{equation}
\end{Proposition}

\noindent{\bf Proof:}
Insert the identity (\ref{up1eq}) inside (\ref{recureq})
and the result follows immediately.
\qed

The case $m=3$ is Dixon's Summation Theorem for terminating
${ }_3F_2$ hypergeometric series.

\begin{Proposition}(Dixon)
$\forall a_1,a_2,a_3\in\mathbb{Z}$,
\[
\ups_3(a_1,a_2,a_3)=
\frac{(2a_1)!(2a_2)!(2a_3)!(a_1+a_2+a_3)!}{(a_1+a_2)!
(a_1+a_3)!(a_2+a_3)!a_1! a_2!a_3!}\ .
\]
\end{Proposition}

\noindent{\bf Proof:}
We insert (\ref{up2eq}) inside (\ref{recureq}) and get
\[
\ups_3(a_1,a_2,a_3)=\sum_l
\frac{(2a_1)!(2a_2)!}{(a_1+a_2)!(2l)!(a_1-l)!(a_2-l)!}
\times
\frac{(2l)!(2a_3)!}{(l+a_3)!l!a_3!}\ .
\]
Since the factor $\frac{1}{l!}$ includes $\bbone\{l\ge 0\}$ one may cancel $(2l)!$ above and below.
Similarly the factors $(2a_2)!$ and $(2a_3)!$ provide the condition $a_2+a_3\ge 0$
which allows us to multiply and divide by $(a_2+a_3)!$. After regrouping we obtain
\[
\ups_3(a_1,a_2,a_3)=
\frac{(2a_1)!(2a_2)!(2a_3)!}{(a_1+a_2)!a_3!a_1!(a_2+a_3)!}
\times\sum_l
\left(\begin{array}{c} a_1 \\ l \end{array}\right)
\left(\begin{array}{c} a_2+a_3 \\ a_2-l \end{array}\right)\ .
\] 
Finally, summing over $l$ using Proposition \ref{ChuV}
and cleaning up the resulting expression gives the desired identity.
\qed

\begin{Remark}
The particular case $a_1=a_2=a_3$ of Dixon's Theorem was conjectured by Morley and proved
in~\cite{Dixon1}. The general case is also due to Dixon~\cite{Dixon2} although it is
sometimes attributed to Fjeldstad~\cite{Fjeldstad}. 
It has also been proved by
Racah~\cite[Appendix A]{Racah}.
There are of course many proofs available, but to the best of our knowledge this is perhaps the first proof which only uses the Chu-Vandermonde identity in both directions, in the spirit of Racah's derivation of his single sum formula for $6j$
symbols~\cite[Appendix B]{Racah}.
\end{Remark}

For $m>3$ there are no more simple closed formulas for $\ups_m(a_1,\ldots,a_m)$,
since, even when $a_1=\cdots=a_m$, no such formula exists as shown by
de Bruijn~\cite[\S4.7]{deBruijn}.
Nevertheless the following result will be enough to prove our main theorem.

\begin{Proposition}
$\forall m\ge 2$, $\forall a_1,\ldots,a_m\ge 0$, we have
$\ups_m(a_1,\ldots,a_m)>0$.
\end{Proposition}

\noindent{\bf Proof:}
For $m=2$ this follows from the explicit formula in (\ref{up2eq}).
The general case follows by a simple induction on $m$.
Namely, for $m\ge 3$, apply (\ref{recureq})
and bound the sum from below by the single term with $l=0$.
This gives
\[
\ups_m(a_1,\ldots,a_m)\ge \frac{(2a_1)!\ (2a_2)!}{(a_1+a_2)a_1!a_2!}
\times\ups_{m-1}(0,a_3,\ldots,a_m)
>0\ .
\]
\qed

We now relate the combinatorial sums of this section to the numerical quantities
$\mathscr{N}_{k,r}$ from \S2.

\begin{Proposition}\label{nonvan}
For all integers $p,q$ with $1\le q\le p$,
\[
\mathscr{N}_{2p,2q+1}=
\frac{(-1)^{p+q}\ (2p)!^{2q+1}}{\prod_{\nu=0}^{2q}[2(p-q+\nu)]!}
\times
\ups_{2q+1}(p-q,p-q+1,\ldots,p+q+1)\ .
\]
\end{Proposition}

\noindent{\bf Proof:}
In the formula (\ref{Ndefeq}) for $\mathscr{N}_{k,r}$ with $k=2p$ and $r=2q+1$ 
and for $0\le j\le k-r+1$,
we write
\[
\left(\begin{array}{c} k \\ j \end{array}\right)
\left(\begin{array}{c} k \\ j+1 \end{array}\right)
\cdots
\left(\begin{array}{c} k \\ j+r-1 \end{array}\right)
=\frac{k!^r}{D}
\]
where the denominator is
\begin{eqnarray*}
D & = & \prod_{\nu=0}^{r-1}\left[
(j+\nu)!(k-j-\nu)!
\right] \\
 & = & \prod_{\nu=0}^{r-1}
(j+\nu)!
\times
\prod_{\nu=0}^{r-1}
(k-j-\nu)! \\
 & = & \prod_{\nu=0}^{r-1}
(j+\nu)!
\times
\prod_{\nu=0}^{r-1}
(k-j-r+1+\nu)!
\end{eqnarray*}
where in the last product we reversed the order of factors, i.e., changed $\nu$
to $r-1-\nu$.
As a result,
\begin{eqnarray*}
\frac{k!^r}{D} & = & \frac{k!^r}{\prod_{\nu=0}^{r-1}
\left[(j+\nu)!
(k-j-r+1+\nu)!\right]} \\
 & = & \frac{k!^r}{\prod_{\nu=0}^{r-1}
(k-r+1+2\nu)!}\times
\prod_{\nu=0}^{r-1}
\left(\begin{array}{c} k-r+1+2\nu \\ j+\nu \end{array}\right)\ .
\end{eqnarray*}
We then change the summation index from $j$ to $i=j+q-p$, write $k,r$ in terms of $p,q$, use the definition of $\ups_m$ and clean up the final expression to get the desired identity.
\qed

The proof of Theorem \ref{mainthm} is now complete.

\section{Explicit examples}

In this section we provide explicit computations for some of our invariants
$\mathscr{P}$. We will assume familiarity with the classical symbolic method and notation, as explained in~\cite[\S2]{Abdesselam12}. 
A superficial glance at the following calculations may give the impression that they are mindless and straightforward. They are not! They in fact provide a rare opportunity for the reader to get ``hands on'' experience with Gordan's 1868 reduction algorithm~\cite{Gordan}.
The latter uses as a guide the graph-theoretic notion of $p$-edge-connectivity.
One of the main ideas is to use variants of the Grassmann-Pl\"{u}cker relation in order to create pairs of vertices with a high number of edges joining them. This means creating factors $(ab)^p$ with $p>\frac{d}{2}$. The highest such exponent in a graph or symbolic bracket monomial is also called the ``grade'' (see~\cite{Weyman} for a modern perspective).
When present, such factors $(ab)^p$ result in the vertices $a$ and $b$ having less than $\frac{d}{2}$
remaining connections to the rest of the graph. Applying the Clebsch-Gordan identity between these remaining edges emanating from $a$ and $b$ only produces transvections
$(\cdot,\cdot)_q$ with $q<\frac{d}{2}$ which is important for the success of the reduction process. The most technically delicate part of Gordan's algorithm is the special treatment of forms of degree divisible by $4$ which have a degree two covariant of the same order as the form itself, namely $(ab)^{\frac{d}{2}}a_x^{\frac{d}{2}}b_x^{\frac{d}{2}}$. In order to eliminate the troublesome factors $(ab)^{\frac{d}{2}}$, Gordan had to consider covariants of degree three and he discovered clever identities~\cite[\S7]{Gordan} which allowed him to improve the grade from $\frac{d}{2}$ to some higher exponent. Our invariants $\mathscr{P}_{k,n}$ featuring in Theorem \ref{mainthm} are entirely made of factors
$(ab)^{\frac{d}{2}}$ and thus take us to the very technical heart of Gordan's method. We do not exactly use his identities but resort to grade-improving identities that are similar in spirit. The most important one is Eq. (\ref{B44eq}) below. Finally, note that the following computations were first done graphically and then converted to the more compact classical symbolic notation. The reader is strongly encouraged to follow the calculations by drawing directed graphs ``on the side'' where symbolic letters $a$, $b$, etc. correspond to vertices and brackets $(ab)$, $(ac)$, etc. correspond to directed edges.

In general, when $d=2k$ with $k$ even and for $p=2,\ldots,k+1$, we have the ``cyclic'' symbolic formula
\[
\mathscr{P}_{k,p}(F)=(a^{(1)}a^{(2)})^k(a^{(2)}a^{(3)})^k(a^{(3)}a^{(4)})^k
\cdots(a^{(p-1)}a^{(p)})^k(a^{(p)}a^{(1)})^k
\]  
where $a^{(1)},\ldots,a^{(p)}$ are symbolic letters for the form $F$.
The particular case of the binary quartic is rather trivial, since we get:
\[
\mathscr{P}_{2,2}=(ab)^2(ba)^2=(F,F)_{4}\ \ ,
\ \ 
\mathscr{P}_{2,3}=(ab)^2(bc)^2(ca)^2=(F,(F,F)_2)_4\ .
\]
It has been known since the middle of the 19th century that these two invariants generate ${\rm Inv}_4$ and are algebraically independent. Thus the HSOP and regular sequence properties hold trivially. Note that the vanishing of $\mathscr{P}_{2,2}$ detects equianharmonic configurations of four points on $\mathbb{P}^1(\mathbb{C})$, while that of the catalecticant
$\mathscr{P}_{2,3}$ detects binary quartics which are limits of sums of two fourth powers of linear forms.

Let us now focus on the more challenging case $d=8$ and explicitly compute the invariants
$\mathscr{P}_{4,2},\mathscr{P}_{4,3},\mathscr{P}_{4,4},\mathscr{P}_{4,5}$ relevant to
our Theorem \ref{mainthm}, in terms of the invariants listed by Shioda~\cite{Shioda}.
In the latter article, Shioda proved that ${\rm Inv}_8$ is generated by elements $J_2,J_3,\ldots,J_{10}$ of respective degree $2,3,\ldots,10$. He showed $J_2,\ldots,J_7$ are algebraically independent and also proved that ${\rm Inv}_8$ is a free module over
$\mathbb{C}[J_2,\ldots,J_7]$ with basis given by $1,J_8,J_9,J_{10},J_9^2$. Namely, he explicitly proved the Cohen-Macaulay property for ${\rm Inv}_8$. From his list, we will need
\begin{eqnarray*}
J_2 & = & (F,F)_8=(ab)^8\ ,\\
J_3 & = & (F,(F,F)_4)_8=(ab)^4(ac)^4(bc)^4\ ,\\
J_4 & = & ((F,F)_6,(F,F)_6)_4\ ,\\
J_5 & = & ((F,(F,F)_6)_4,(F,F)_6)_4=(ab)^6(ac)^2(bc)^2(cd)^2(ce)^2(de)^6\ .
\end{eqnarray*}
These same invariants are respectively denoted by $A,B,C, f_{k,k}$ in~\cite{vonGall1}.
We clearly have $\mathscr{P}_{4,2}=J_2$ and $\mathscr{P}_{4,3}=J_3$. However, the computation of $\mathscr{P}_{4,4}$ and $\mathscr{P}_{4,5}$ is more involved and necessitates a detailed study of degree three covariants of order four and eight.

For $\lambda=(\lambda_1,\lambda_2,\lambda_3)$, with $\lambda_1\ge\lambda_2\ge\lambda_3\ge0$, 
an integer partition of $4$, we define
\[
A_{\lambda}=(ab)^{\lambda_3+2}(ac)^{\lambda_2+2}(bc)^{\lambda_1+2}
a_x^{\lambda_1}b_x^{\lambda_2}c_x^{\lambda_3}\ .
\]
It is easy to that the $A_{\lambda}$ linearly generate the space of covariants of degree $3$ and order $4$. Note that we have the basic identity
\[
a_x(bc)=b_x(ac)-c_x(ab)
\]
which allows one to rewrite any of the above three terms using the other two.
Consider
$A_{31}=(ab)^{2}(ac)^{3}(bc)^{5}a_x^3 b_x$
and use the basic identity to rewrite $b_x(ac)$. We then get
\[
A_{31}=(ab)^{2}(ac)^{2}(bc)^{6}a_x^4+(ab)^{3}(ac)^{2}(bc)^{5}a_x^3 c_x\ .
\]
Exchanging the symbolic letters $b$ and $c$ in the second term and using the antisymmetry of brackets $(cb)=-(bc)$, one obtains the relation $A_{31}=A_4-A_{31}$ and thus
$A_{31}=\frac{1}{2} A_4$.
Rewriting $a_x(bc)$ inside
$A_{211}=(ab)^{3}(ac)^{3}(bc)^{4}a_x^2 b_x c_x$ gives the relation $A_{211}=-2A_{211}$
and thus $A_{211}=0$.
Finally, rewriting $a_x(bc)$ inside $A_{22}=(ab)^{2}(ac)^{4}(bc)^{4}a_x^2 b_x^2$ gives
\begin{equation}
A_{22}=A_{31}+A_{211}=\frac{1}{2}A_4=\frac{1}{2}(F,(F,F)_6)_4\ .
\label{A22eq}
\end{equation}
We now consider partitions $\lambda=(\lambda_1,\lambda_2,\lambda_3)$ of $8$ and define
\[
B_{\lambda}=(ab)^{\lambda_3}(ac)^{\lambda_2}(bc)^{\lambda_1}
a_x^{\lambda_1}b_x^{\lambda_2}c_x^{\lambda_3}\ .
\]
These bracket monomials linearly generate the space of
of covariants of degree $3$ and order $8$.
Rewriting $c_x(ab)$ inside $B_{521}=(ab)(ac)^{2}(bc)^{5}a_x^5 b_x^2 c_x$
gives
$B_{521}=B_{53}-B_{62}$.
Rewriting $a_x(bc)$ inside $B_{422}=(ab)^2(ac)^{2}(bc)^{4}a_x^4 b_x^2 c_x^2$
gives
$B_{422}=2B_{332}$.
Rewriting $a_x(bc)$ inside $B_{431}=(ab)(ac)^{3}(bc)^{4}a_x^4 b_x^3 c_x$
gives $B_{431}=-B_{431}-B_{332}$ and thus
$B_{332}=-2B_{431}$.
However, one can also rewrite $b_x(ac)$ inside $B_{431}$ again which results in the relation
\[
B_{431}=B_{521}+B_{422}=B_{521}-4B_{431}
\]
and therefore
\[
B_{431}=\frac{1}{5}B_{521}=\frac{1}{5}B_{53}-\frac{1}{5}B_{62}\ .
\]
Finally, rewriting $a_x(bc)$ inside $B_{44}=(ac)^{4}(bc)^{4}a_x^4 b_x^4$
gives $B_{44}=B_{53}+B_{431}$, i.e.,
\begin{equation}
B_{44}=\frac{6}{5}B_{53}-\frac{1}{5}B_{62}\ .
\label{B44eq}
\end{equation}
Now, by definition,
\[
\mathscr{P}_{4,4}=(ab)^4(bc)^4(cd)^4(da)^4=(F,B_{44})_8
=\frac{6}{5}(F,B_{53})_8-\frac{1}{5}(F,B_{62})_8\ .
\]
The last transvectant, which is easier to compute, is given by
\[
(F,B_{62})_8=(\ d_x^8\ ,\ (ac)^{2}(bc)^{6}a_x^6 b_x^2\ )_8
=(da)^{6}(db)^{2}(ac)^{2}(bc)^{6}\ .
\]
We now apply the Clebsch-Gordan (CG) identity~\cite[Eq. 2.9]{Abdesselam12}
to the parallel groups of two brackets each $(db)^2$ and $(ac)^2$.
This results in the formula
\[
(F,B_{62})_8=\sum_{j=0}^{2}\frac{{\binom{2}{j}}^2}{\binom{5-j}{j}}
((F,F)_{6+j},(F,F)_{6+j})_{4-2j}\ .
\]
Since odd transvectants of a form with itself vanish, only the terms $j=0,2$ survive and the resulting expression readily simplifies to
\[
(F,B_{62})_8=J_4+\frac{1}{3}J_{2}^{2}\ .
\]
Similarly,
$(F,B_{53})_8=(da)^{5}(db)^{3}(ac)^{3}(bc)^{5}$
to which we apply the CG identity for the bracket groups $(db)^3$ and $(ac)^3$ with the result
\[
(F,B_{53})_8=\sum_{j=0}^{3}\frac{{\binom{3}{j}}^2}{\binom{7-j}{j}}
((F,F)_{5+j},(F,F)_{5+j})_{6-2j}\ .
\]
Only the $j=1,3$ terms survive and the sum simplifies to
$(F,B_{53})_8=\frac{3}{2}J_4+\frac{1}{4}J_{2}^{2}$.
Hence
\[
\mathscr{P}_{4,4}=\frac{8}{5}J_4+\frac{7}{30}J_{2}^{2}\ .
\]
We now consider
$\mathscr{P}_{4,5}=(ab)^4(bc)^4(cd)^4(de)^4(ea)^4$
and apply the CG identity to the bracket groups $(ea)^4$ and $(cd)^4$. This gives
\begin{eqnarray}
\mathscr{P}_{4,5} & = & \sum_{j=0}^{4}\frac{{\binom{4}{j}}^2}{\binom{9-j}{j}}
\left(\ (F,F)_{4+j}\ , (ab)^{4}(bc)^{4}(ac)^{j}a_x^{4-j}c_x^{4-j}\ \right)_{8-2j}
\nonumber \\
 & = & ((F,F)_4,B_{44})_8+\frac{12}{7}((F,F)_6,A_{22})_4+\frac{1}{5}J_2 J_3\ .
\label{P45breakeq}
\end{eqnarray}
The middle term in (\ref{P45breakeq}) is easy to handle thanks to (\ref{A22eq}). Indeed,
\[
((F,F)_6,A_{22})_4=\frac{1}{2}((F,F)_6,A_{4})_4=\frac{1}{2}J_5\ .
\]
Using (\ref{B44eq}), the computation of the first term in (\ref{P45breakeq})
reduces to that of $((F,F)_4,B_{53})_8$ and $((F,F)_4,B_{62})_8$.
We have
\[
((F,F)_4,B_{62})_8=\left(\ (de)^4 d_x^4 e_x^4\ ,\ (ac)^2(bc)^6a_x^6b_x^2\ \right)_8
=\frac{2\times 4\times 3\times 6!}{8!}\ C_1
+\frac{8\times 4\times 6!}{8!}\ C_2
\]
where
\[
C_1=
(de)^4(da)^2(db)^2(ea)^4(ac)^2(bc)^6\ \ ,
\ \ C_2=
(de)^4(da)^3(db)(eb)(ea)^3(ac)^2(bc)^6\ .
\]
This formula is obtained by expanding the transvectant into bracket monomials,
similarly to summing over Wick contractions in quantum field theory followed by collecting ``topologically equivalent graphs''. Here, this means one has to symmetrize over the matchings of the $8$ factors carrying an $x$ on the left to the $8$ ones on the right. The combinatorial weights can be computed ``probabilistically'' as follows. The two $b_x$ factors can either connect to the same $d$ or $e$ vertex or to each of the two vertices. The $C_1$ and $C_2$ terms above respectively correspond to these two possibilities. For the $C_1$ term one has
$2$ choices for which vertex receives the two connections. Say it is $d$. Then one has $4$ choices to connect the first $b_x$ and $3$ more choices for the second one. Finally there are $6$ remaining $x$-carrying factors on the left which will all
connect to the $a$ vertex and this accounts for $6!$ more choices. The numerator for the $C_2$ term is obtained by a similar counting argument.

By applying the CG identity to the bracket groups $(db)^2$ and $(ac)^2$ one finds
\begin{eqnarray*}
((F,F)_4,C_1)_8 & = & \sum_{j=0}^{2}\frac{{\binom{2}{j}}^2}{\binom{5-j}{j}}
\left(\ 
(de)^4(ea)^4(da)^{2+j}d_x^{2-j}a_x^{2-j}
\ ,\ 
(bc)^{6+j}b_x^{2-j}c_x^{2-j}
\ \right)_{4-2j} \\
 & = & (A_{22},(F,F)_6)_4+\frac{1}{3}J_3 J_2=\frac{1}{2}J_5+\frac{1}{3}J_2 J_3
\end{eqnarray*}
since the $j=1$ term vanishes.

The $C_2$ term needs more care since we will now apply the CG identity to the pair of bracket groups $(db)(eb)$ and $(ac)^2$. This gives, in redundant but hopefully self-explanatory form,
\[
C_2=
\frac{{\binom{2}{0}}^2}{\binom{5-0}{0}}
(A_{211},(F,F)_6)_4+
\frac{{\binom{2}{2}}^2}{\binom{5-2}{2}}
(J_3,(F,F)_8)_0
\]
since, again the $j=1$ term is absent because $(F,F)_7=0$.
Since we have shown that $A_{211}=0$, we get $C_2=\frac{1}{3}J_2 J_3$ and therefore
\[
((F,F)_4,B_{62})_8=\frac{3}{7}\left(\frac{1}{2}J_5+\frac{1}{3}J_2 J_3\right)
+\frac{4}{7}\left(\frac{1}{3}J_2 J_3\right)=\frac{3}{14}J_5+\frac{1}{3}J_2 J_3\ .
\] 

The ``probabilistic'' expansion of the transvectant into bracket monomials gives
\[
((F,F)_4,B_{53})_8=\frac{1}{7}\ C_3+\frac{6}{7}\ C_4
\]
with
\[
C_3=(de)^4(da)(db)^3(ea)^4(ac)^3(bc)^5\ \ ,
\ \ 
C_4=(de)^4(da)^2(db)^2(ea)^3(eb)(ac)^3(bc)^5\ .
\]
Applying the CG identity to $(db)^3$ and $(ac)^3$ in $C_1$ gives
\[
C_3=\frac{{\binom{3}{1}}^2}{\binom{6}{1}}
(A_{22},(F,F)_6)_4+\frac{{\binom{3}{3}}^2}{\binom{4}{3}}
J_2 J_3=\frac{3}{4}J_5+\frac{1}{4}J_3 J_2\ .
\]
Applying the CG identity to $(db)^2(eb)$ and $(ac)^3$ in $C_4$, with extra care as we did for $C_2$,
we obtain
\[
C_4=\frac{{\binom{3}{1}}^2}{\binom{6}{1}}
\left(\ \frac{1}{3}A_{22}+\frac{2}{3}A_{211}\ ,\ (F,F)_6\ \right)_4
+\frac{{\binom{3}{3}}^2}{\binom{4}{3}} J_3 J_2
=\frac{1}{4}J_5+\frac{1}{4}J_2 J_3\ .
\]
Substituting back, we get
\[
((F,F)_4,B_{53})_8=\frac{9}{28}J_5+\frac{1}{4}J_2 J_3\ \ ,
\ \ 
((F,F)_4,B_{44})_8=\frac{12}{35}J_5+\frac{7}{30}J_2J_3\ ,
\]
and finally
\[
\mathscr{P}_{4,5}=\frac{6}{5}J_5+\frac{13}{30}J_2 J_3\ .
\]
We thus explicitly showed that $J_2,J_3,J_4,J_5$ and
$\mathscr{P}_{4,2},\mathscr{P}_{4,3},\mathscr{P}_{4,4},\mathscr{P}_{4,5}$
are triangularly related and established the following result.

\begin{Proposition}\label{octavicprop}
Conjecture \ref{genconj} holds for $(d,n)=(4,2)$ and $(8,4)$.
\end{Proposition}

\begin{Remark}
To the best of our knowledge, there are only two proofs of Dixmier's conjecture for the octavic. In~\cite[p. 137]{Dixmier}, Dixmier says that the $d=8$ case was checked by computer by Bartels and then he later also checked it ``by hand''.
However, no details are given. The second proof~\cite{LittelmannP} is also computer-assisted.
Incidentally, one can use the argument of~\cite{LittelmannP} to quickly derive Proposition \ref{octavicprop} from Theorem \ref{mainthm}. From the knowledge of the first few terms of the Hilbert series, and the algebraic independence of
$\mathscr{P}_{4,2},\mathscr{P}_{4,3},\mathscr{P}_{4,4},\mathscr{P}_{4,5}$ one deduces that they are triangularly related to $J_2,J_3,J_4,J_5$ and therefore constitute a regular sequence in ${\rm Inv}_8$. Of course, this rests crucially on the prior knowledge of a HSOP
with the right degrees which here is due to Shioda~\cite{Shioda}.
\end{Remark}

\section{Why Dixmier's conjecture is not only hard but ridiculously so}

We again assume the degree of our generic binary form is $d=2k$ with $k$ even.
In order to understand the invariants $\mathscr{P}_{n,p}$ it is important
to be able to see them.
This is provided by the Feynman diagram calculus developed in~\cite{Abdesselam12}, with additional explanations given in~\cite[\S4.2]{AbdesselamC18}. We will not repeat here the corresponding definitions and refer the reader to these two resources.
In the notations of~\cite[\S2]{Abdesselam12}
one has
\[
\mathscr{P}_{n,p}(F)= 
\parbox{2cm}{
\psfrag{1}{$\scriptstyle{k}$}
\psfrag{2}{$\scriptstyle{n}$}
\psfrag{3}{$\scriptstyle{n-k}$}
\psfrag{4}{$F$}
\raisebox{1.7ex}{\includegraphics[width=2cm]{Fig1.eps}}}
\qquad\qquad\qquad\qquad\qquad\qquad
=
\parbox{2cm}{
\psfrag{1}{$\scriptstyle{2k}$}
\psfrag{2}{$\scriptstyle{n}$}
\psfrag{4}{$F$}
\includegraphics[width=2cm]{Fig2.eps}}
\qquad\qquad\qquad\qquad\qquad
\]

\bigskip\noindent
The first graphical formula uses the ``microscopic'' notation of~\cite[\S2]{Abdesselam12}
whereas the second one uses
the more compact ``macroscopic'' notation of the same reference.
There is no numerical ambiguity (overall sign or proportionality factor) in both
graphical formulas.

In the simplest quadratic case one has, for all $n\ge k$, the formula
\[
\mathscr{P}_{n,2}(F)=
\frac{(n-k)!\ (n+k+1)!\ k!^2}{n!^2\ (2k+1)!}\times (F,F)_{2k}\ .
\]
The latter is an immediate consequence of the graphical form of an identity of
Clebsch~\cite[Eq. (2.10)]{Abdesselam12}.
However, for the cubic invariant case
one has that $\mathscr{P}_{n,3}(F)=\alpha_{n,k}\beta_{n,k} \mathscr{P}_{k,3}(F)$ where
$\mathscr{P}_{k,3}(F)=(F,(F,F)_k)_{2k}$, $\alpha_{n,k}$ is a trivial nonvanishing factor ($\pm$ a square root of a ratio of products of factorials) and $\beta_{n,k}$ is the Wigner $6j$ symbol
\[
\beta_{n,k}=
\left\{
\begin{array}{ccc}
k & k & k \\
\frac{n}{2} & \frac{n}{2} & \frac{n}{2}
\end{array}
\right\}\ .
\]
This computation is an example of star-triangle relation in the quantum theory of angular momentum. It follows from similar manipulations as
in~\cite[Proof of Theorem 7.2]{AbdesselamC12}. The definition of $6j$ symbols in relation to the invariant theory of binary forms is recalled in~\cite[\S7]{AbdesselamC9}.
A trivial consequence of Conjecture \ref{genconj} is that $\beta_{n,k}$ should not vanish.
Encouraged by computer checks due to J. Van der Jeugt, we propose the following purely combinatorial conjecture.

\begin{Conjecture}\label{6jconj}

For all pairs of integers $(k,n)$ with $n\ge k\ge 2$ and with the exception of $(k,n)=(2,3)$,
one has
\[
\left\{
\begin{array}{ccc}
k & k & k \\
\frac{n}{2} & \frac{n}{2} & \frac{n}{2}
\end{array}
\right\}\neq 0\ .
\] 
\end{Conjecture}

Note that we did not assume $k$ even in the last statement. Simply by requiring the nonvanishing of the degree three invariant in Dixmier's list, it is immediate that
Dixmier's conjecture would imply the particular cases $n=2k-2$ with $k$ even and $k\ge 2$
of Conjecture \ref{6jconj}.
By Racah's
formula for $6j$ symbols~\cite[Appendix B]{Racah} (see also~\cite[\S7]{AbdesselamC12}), 
Conjecture \ref{6jconj}
amounts to the nonvanishing of the combinatorial sum
\[
\sum_{j}(-1)^j
\left(\begin{array}{c} j+1 \\ 3k+1 \end{array}\right)
\left(\begin{array}{c} k \\ j-k-n \end{array}\right)^3
\]
for the same range of pairs $(k,n)$.
A matrix plot of the Mathematica command
\begin{center}
\texttt{Table[Sign[SixJSymbol[\{2 + i, 2 + i, 2 + i\}, \{(2 + i + j)/2, (2 + i + j)/2,}\\
\texttt{(2 + i + j)/2\}]], \{i, 0, 200\}, \{j, 0, 200\}]}
\end{center}
produces the following picture which exhibits a rather intricate sign pattern.
\[
\parbox{14cm}{
\includegraphics[height=14cm]{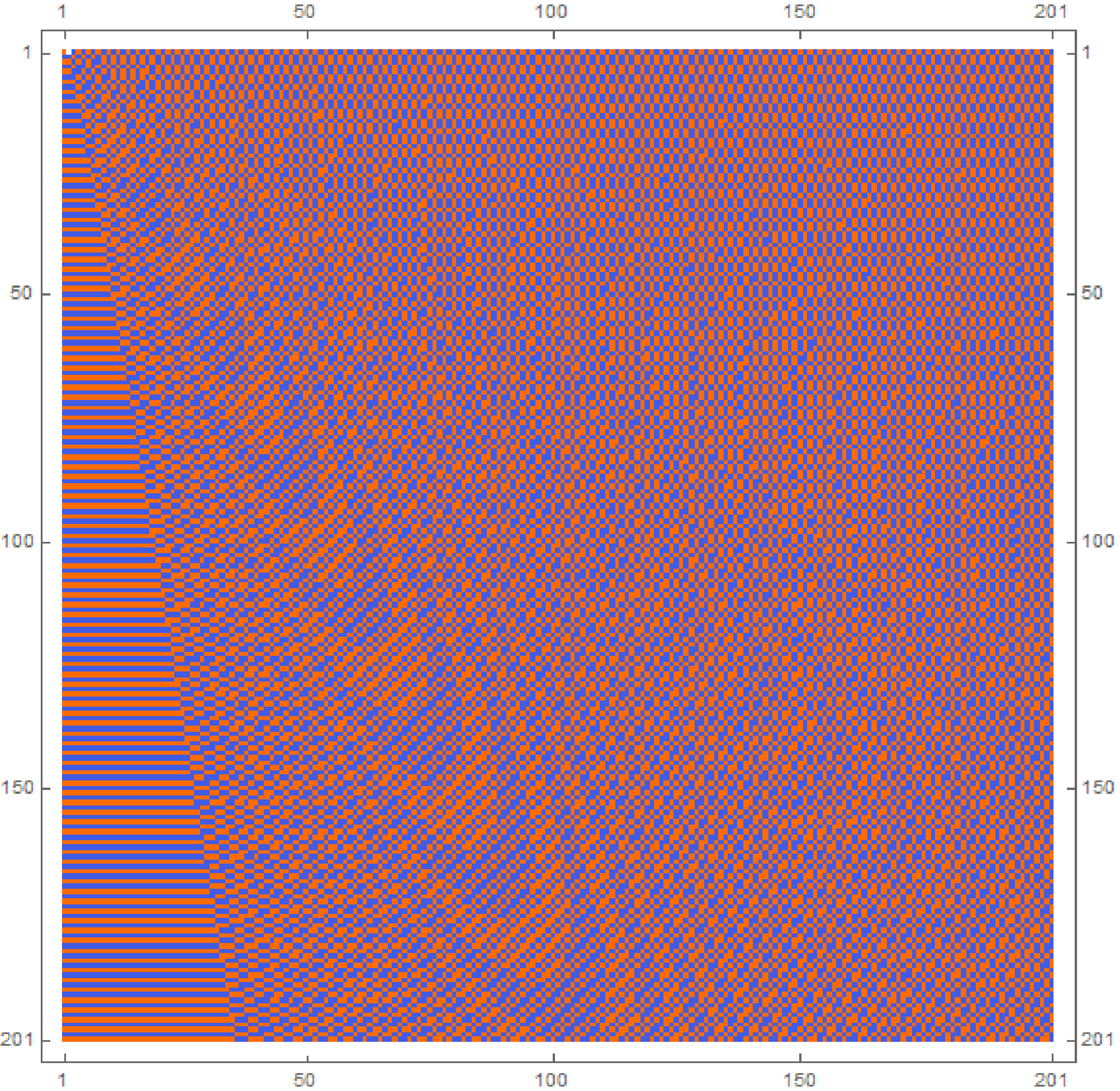}
}
\]
If $r=1,\ldots,201$ is the row number from top to bottom and if $c=1,\ldots,201$ is the column number from left to right then the coloring of the square in $r$-th row and $c$-th
column of the previous picture gives the sign of the $6j$ symbol
\[
\left\{
\begin{array}{ccc}
r+1 & r+1 & r+1 \\
\frac{r+c}{2} & \frac{r+c}{2} & \frac{r+c}{2}
\end{array}
\right\}\ .
\]
Light square means positive, dark means negative and white square means zero.
In general, zeros of $6j$ symbols are poorly understood (see,
e.g.,~\cite{SrinivasaRaoRK,VanderJeugt,CaglieroS}). As to the curious $(k,n)=(2,3)$
exception or white square in the top left corner of the picture, it has been given a representation-theoretic explanation in~\cite{Minnaert}.
Needless to say, the author of the present article has no idea how to prove
Conjecture \ref{6jconj}, even in the Dixmier case where $n=2k-2$.
We also learned from C. Krattenthaler that for simpler-looking combinatorial sums, similar nonvanishing conjectures are open
(see, e.g.,~\cite[Conjecture A]{CarnevaleV}).

\bigskip
\noindent{\bf Acknowledgements:}
{\small
For useful discussions or correspondence the author thanks A. Brouwer, J. Chipalkatti, C. Huneke, C. Krattenthaler and
J. Van der Jeugt. The author also thanks the anonymous referee for suggesting useful improvements.
}

\medskip
\noindent{\bf Conflict of interest statement:}
{\small
The author states that there is no conflict of interest.
}

\end{document}